\documentclass[12pt,reqno]{amsart}

\usepackage{euscript}
\usepackage{amsmath}
\usepackage{epsfig}
\usepackage{amssymb}
\usepackage{epic}

\theoremstyle{plain}

\theoremstyle{definition}

\newcommand{\bc}{{\mathbb C}}
\newcommand{\bp}{{\mathbb P}}
\newcommand{\Z}{{\mathbb Z}}
\newcommand{\bn}{{\mathbb N}}
\newcommand{\bq}{{\mathbb Q}}
\newcommand{\vp}{\varphi}

\newcommand{\bH}{{\mathbb H}}

\newcommand{\calO}{{\mathcal{O}}}
\newcommand{\calC}{{\mathcal{C}}}
\newcommand{\calS}{{\mathcal{S}}}
\newcommand{\calt}{{\mathcal T}}

\newcommand{\cc}{\bar{C}}
\newcommand{\ssw}{{\bf sw}}
\newcommand{\et}{\EuScript{T}}
\newcommand{\Gammma}{{G}}
\newcommand{\calv}{{\mathcal V}}
\newcommand{\cale}{{\mathcal E}}

\newcommand{\Q}{{\mathbb Q}}
\newcommand{\Hom}{{\rm Hom}}
\def\R{\mathbb R}

\begin{document}

\title{On rational cuspidal curves, open surfaces and local singularities}
\author{J. Fern\'andez de Bobadilla,
I. Luengo, A. Melle-Hern\'andez, \and A. N\'emethi}
\address{Departamento de Matem\'aticas Fundamentales\\
\\ Facultad de Ciencias U.N.E.D.\\ c/ Senda del Rey 9, 28040 Madrid, Spain.
}
\address{Facultad de Matem\'aticas\\ Universidad Complutense\\
Plaza de Ciencias\\ E-28040, Madrid, Spain}
\address{Department of Mathematics\\Ohio State University\\Columbus,
OH 43210,USA; and R\'enyi Institute of Mathematics, Budapest, Hungary.}

\email{javier@mat.uned.es}
\email{iluengo@mat.ucm.es}
\email{amelle@mat.ucm.es}
\email{nemethi@math.ohio-state.edu; nemethi@renyi.hu}

\thanks{\noindent The first author is partially supported by Ram\'on y Cajal contact.
The first three authors thank  Marie Curie
Fellowship for the Transfer of Knowledge supporting their visit
at the R\'enyi Institute, Budapest, Hungary; the first three authors 
are partially supported by MTM2004-08080-C02-01;
the last author is also supported by the Marie Curie
Fellowship  and OTKA Grant.}




\maketitle
\pagestyle{myheadings}
\markboth{{\normalsize J. Fern\'andez de Bobadilla, I. Luengo, A. Melle-Hern\'andez, A.
N\'emethi}}{
{\normalsize Rational cuspidal curves, open surfaces and singularities}}

{\small

\section{Introduction}

Let $C$ be an irreducible  projective  plane curve
in the complex projective space $\bp^2$. The classification
of such curves, up to the action of the automorphism group
$PGL(3,\bc)$ on $\bp^2$,
is a very difficult  open problem with many interesting connections.
The main goal is to determine, for a given $d$, whether
there exists a projective plane curve of degree $d$ having a fixed
number of singularities of given topological type.
In this note we are mainly  interested in the case when $C$ is a
rational curve.

The problem remains  very difficult even if we aim much less, 
e.g. the determination of the maximal number of cusps among all the 
rational cuspidal plane curves
(a problem proposed by F. Sakai in \cite{open}),
 ---\,  this number  is expected to be small. In 
\cite{Tono2}  K.~Tono  recently proved  that it is less than 9; the maximal 
number of cusps known by the authors is 4; and, in fact, it is expected to
be 4. The referee pointed out to us that in mid-90-s S.~Orevkov found a bound bigger than 8 but he never published his result. 

This remarkable  problem of classification is not only important for its own 
sake, but it is also
connected with crucial properties, problems and conjectures in the theory of
open surfaces, and in the classical algebraic geometry.

For instance, the open surface
$\bp^2\setminus C$ is $\bq$-acyclic if and only if $C$ is a rational
cuspidal curve. On the other hand, regarding these surfaces,
Flenner and Zaidenberg in \cite{FZ0}  formulated  the
\emph{rigidity conjecture}. This says that
every $\bq$-acyclic affine surfaces $Y$
with logarithmic Kodaira dimension $\bar{\kappa}(Y)=2$
must be rigid. This conjecture for $Y=\bp^2\setminus C$ would imply
the projective rigidity of  the curve $C$
in the sense that every equisingular deformation of $C$ in $\bp^2$ would be
projectively equivalent to $C$.
(Notice that  if $C$ has at least three cusps
then  $\bar{\kappa}(\bp^2\setminus C)=2$ by \cite{Wak};
and, in fact, all curves with  $\bar{\kappa}(\bp^2\setminus C)<2$
are classified, see below.) 
 Many known examples support
the rigidity conjecture, see \cite{Fenske1,FZ0,FZ1,FZ2}.
Zaidenberg in \cite{open} also conjectured that the set of shapes of the Eisenbud-Neumann splice diagrams is finite for all  $\bq$-acyclic affine surfaces $Y$
with $\bar{\kappa}(Y)=2$ (this is stronger than the existence of a uniform bound for the number of cusps for rational cuspidal curves).

Another related, very old, famous open problem has its roots in early 
algebraic geometry, and wears the name of 
Coolidge and Nagata, see \cite{Coo,Na}. It predicts that   
every rational cuspidal curve can be transformed by
a Cremona transformation into a line. 

The aim of this article is to present  some of these
conjectures and related problems, and to complete  them with 
some results and new conjectures from the recent work of the authors.

In section 3 we present the 
Nagata-Coolidge problem. The main theme of section 4 is 
Orevkov's conjecture
\cite{Or}, which formulates an inequality involving the degree 
$d$ and numerical invariants of local singularities. In a different 
formulation, this is equivalent with the positivity of the virtual
dimension of the space of curves with fixed degree and certain local type 
of singularities which can be geometrically realized.
(This was used, as a `first test',  by the second author
to check that some singularities might be realized or not.) 
 The equivalence of the 
two inequalities is proved via some  properties of the  numerical, local and 
global deformation invariants; this material is presented in 
sections 2 and 3. 

Section 5 starts with some classification results:
the classification of projective plane curves with 
$\bar{\kappa}(\bp^2\setminus C)<2$ by the work of Kashiwara, Lin, Miyanishi, 
Sugie,  Tsunoda, Tono, Wakabayashi, Yoshihara,  Zaidenberg (among others).
Also, we recall the classification (of the authors)
of rational unicuspidal curves whose cusp has one Puiseux pair. The end of
this section deals with the {\em rigidity conjecture} of Flenner and
Zaidenberg. 

In section 6 we present  the author's `compatibility property'
(a sequence of inequalities), conjecturally
satisfied by the degree and local invariants of the singularities   of a
rational cuspidal curve. It turns out that in the unicuspidal case,  the 
inequalities are true if and only if they are (in fact) equalities. 
One of the reformulations (valid in the unicuspidal case)
of this conjectural series of identities is the
 {\em semigroup distribution property}, which is a very precise 
compatibility property connecting the semigroup of the local singularity 
and the degree of the curve. We also present 
one of its equivalent statements 
 suggested to us by Campillo.  Section 7 explains the relation of
the semigroup distribution property (for the unicuspidal case) 
with the `Seiberg-Witten invariant conjecture' (of the forth author and
Nicolaescu \cite{[51]}), which basically leads us to this compatibility 
property. The last section present another connection with the 
Seiberg-Witten theory (based on the articles \cite{[61],[63],[65]}), 
but now exploiting the relation with the 
Heegaard-Floer homology (introduced and studied by Ozsv\'ath and Szab\'o 
\cite{OSz}). 
Here, a  crucial intermediate  object is the `graded root' 
(introduced by the forth author \cite{[61]}), 
which provides a completely different
 interpretation of the semigroup distribution property. 
 
\vspace{2mm}

 The authors thank A. Campillo,  K. Tono and M.G. Zaidenberg for helpful
information, comments and interesting discussions. The authors thank to the referee
for her/his very useful remarks and comments. 

\vspace{1mm}

The first three authors are  grateful for the warm 
hospitality of the Alfr\'ed 
 R\'enyi Institute of Mathematics (Budapest, Hungary) where they spent
a fruitful period in an ideal mathematical environment.

\section{Local data}

 One of our main goals is to characterise the local embedded
topological types of local singular germs $(C,p_i)\subset
(\bp^2,p_i) $ which can be realized by a projective plane curve
$C$ of degree $d$. (Here, the points $\{p_i\}$ are not fixed.)
In this section we collect some results about the (deformation of) 
local invariants $(C,p_i)$.

\subsection{}
Let $(C,p)$ be the germ at $p$ of a reduced curve $C\subset \bp^2$ and
 let $f\in \calO_{\bp^2,p}$ be a function
defining the singularity $(C,p)$ in some local coordinates $x$ and $y$.

\subsubsection{}\label{1.1} Let $\phi:\calC_{(C,p)}\to \calS_{(C,p)}$ be
a \emph{semi-universal deformation} of $(C,p)$. The base space
$\calS_{(C,p)}$ is smooth and its tangent space is isomorphic to
the vector space $\calO_{C,p}/(f_x,f_y).$ In particular its
dimension is equal to the \emph{Tjurina number}
$\tau(C,p):=\dim_\bc(\calO_{\bp^2,p}/I^{ea}(C,p))$, where
$I^{ea}(C,p)$ is the ideal  $(f,f_x,f_y)$ (see~\cite{Lo}.)

\subsubsection{}\label{1.2} Let $\calS_{(C,p)}^{es}\subset \calS_{(C,p)}$ 
be the smooth
subgerm of a \emph{semi-universal equisingular deformation} of $(C,p)$,
see \cite{wahl}, Theorem 7.4. Let $T_\varepsilon:=
Spec(\bc[\varepsilon]/(\varepsilon^2))$
be the base space of \emph{first order infinitesimal} deformations.
The \emph{equisingularity ideal} of $(C,p)$ is the ideal
\begin{equation*}
I^{es}(C,p)=\{g\in\calO_{\bp^2,p}\quad
|\quad  f+\varepsilon g \ \text{is equisingular over } T_\varepsilon \}.
\end{equation*}
It contains $I^{ea}(C,p)$. Moreover the tangent space of $\calS_{(C,p)}^{es}
\subset \calS_{(C,p)}$ is isomorphic to the vector space
$I^{es}(C,p)/I^{ea}(C,p).$ In particular the codimension of
$\calS_{(C,p)}^{es}$ in $\calS_{(C,p)}$ is the topological invariant
\begin{equation*}
\tau^{es}(C,p):=\dim_\bc (\calO_{\bp^2,p}/I^{es}(C,p)). \tag{1}
\end{equation*}
The equisingular stratum $\calS_{(C,p)}^{es}$ coincides with the
$\mu$-\emph{constant
stratum} of $\calS_{(C,p)}$, where
$\mu=\mu(C,p)$ is the Milnor number of $(C,p)$.

\subsubsection{}\label{1.3} Let $n:\calO_{C,p}\hookrightarrow
\tilde{\calO}_{C,p}:=\prod_{i=1}^r\bc\{t\}$ be the normalisation map,
where $r$ is the number of local irreducible components of $(C,p).$
Let $\delta(C,p)$ be the dimension of the cokernel
$\tilde{\calO}_{C,p}/n(\calO_{C,p})$.
Let $\text{cond}(\calO_{C,p})$ be the conductor ideal of $\calO_{C,p}$, that is
the annihilator of the $\calO_{C,p}$-module
$\tilde{\calO}_{C,p}/n(\calO_{C,p})$.

The $\delta$-\emph{constant stratum}
$\calS_{(C,p)}^{\delta}\subset \calS_{(C,p)}$ is the locus of points
where $\delta$ is constant, see
\cite{dh,teissier}. In general, $\calS_{(C,p)}^{\delta}$ is not longer smooth
at $C$,  but it is locally irreducible. Even more, the tangent cone
of  $\calS_{(C,p)}^{\delta}$ at $C$  is
always a linear space which is identified with the vector space
$\text{cond}(\calO_{C,p})/I^{ea}(C,p)$, \cite{dh}, Theorem 4.15.
In particular, the codimension of $\calS_{(C,p)}^{\delta}$ in $\calS_{(C,p)}$ is equal to
\begin{equation*}
\text{codim}\, (\calS_{(C,p)}^{\delta}\subset \calS_{(C,p)})=\delta(C,p).\tag{2}
\end{equation*}

\subsubsection{}\label{1.4} S.~Diaz and J.Harris \cite{dh} showed the inclusion of ideals $$I^{ea}(C,p)\subset I^{es}(C,p)
\subset \text{cond}(\calO_{C,p}) \subset \calO_{C,p}.$$
Since the equisingular stratum $\calS_{(C,p)}^{es}$ is contained in the $\delta$-constant
stratum $\calS_{(C,p)}^{\delta}$ then
\begin{equation*}
\bar{M}(C,p):=\text{codim}\,(\calS_{(C,p)}^{es}\subset \calS_{(C,p)}^{\delta})=\tau^{es}(C,p)-\delta(C,p) \geq 0.\tag{3}
\end{equation*}
In fact, $\bar{M}(C,p)$  is a topological invariant of the singularity
(see next paragraph).
S.~Orevkov and M.~Zaidenberg used $\bar{M}(C,p)$
in a different situation, cf. \cite{Or,orza}.

\subsection{}\label{1.5}  The minimal embedded resolution
of $(C,p)$ is obtained via a sequence of blowing-ups at infinitely
near points to $p$. The
topological numerical invariants $\delta(C,p),\mu(C,p),\tau^{es}(C,p)$
and $\bar{M}(C,p)$ can be written in terms of the
multiplicity sequence    $[m_1^p,\ldots,m_{k_p}^p]$
of $(C,p)$ associated with this minimal resolution.
This is the set of multiplicities
of the strict transform of $C$ at all the infinitely near points
along all the branches. Let $r_p$ denote the number of branches 
of $C$ at $p.$
(We do not omit the $1$'s; 
for this notation we follow \cite{MS} and \cite{FZ0}.)

Using Milnor's formula $\mu(C,p)=2\delta(C,p)-r_p+1$ 
(see e.g. Theorem 6.5.9 in \cite{Wall}) one gets
\begin{equation*}
\mu(C,p)+r_p-1=2\delta(C,p)=\sum_{i=1}^{k_p} m_i^p(m_i^p-1). \tag{4}
\end{equation*}

A blow-up in the minimal good resolution is called \emph{inner} (or
\emph{subdivisional}) if its center is at the intersection point of
two exceptional curves of the resolution process (such a center is
called a \emph{satellite point} in \cite{W}). If the center
is situated on exactly one exceptional divisor, then the blow-up
is called \emph{outer} (or \emph{sprouting}).  Notice that
the first blow-up is neither inner nor outer.
A center is called a \emph{free infinitely near point}
if it is either outer or it is $p$, the center of the very first blow up.
It is convenient to count this
first infinitely near point by two.
 Let $\omega_p$, resp. $\rho_p$,
denote the number of inner,  respectively of  outer, blow-ups.
Then $k_p-1=\omega_p+\rho_p$ and the number of free infinitely near
points is $L_p:=2+\rho_p=k_p-\omega_p+1.$

\medskip

\noindent
C.T.C.~Wall in Theorem 8.1 of  \cite{W} 
 proved the following formula for $\tau^{es}(C,p)$ 
(see also Proposition 11.5.8 in \cite{Wall})
\begin{equation*}
\tau^{es}(C,p)=\sum_{i=1}^{k_p} \frac{(m_i^p-1)(m_i^p+2)}{2}+\omega_p-1=
\sum_{i=1}^{k_p} \frac{m_i^p(m_i^p+1)}{2}- L_p. \tag{5}
\end{equation*}
Different proofs of this formula have been also given by J.F.~Mattei \cite{Mat},  E.~Shustin \cite{Shu} or Theo de Jong \cite{dJ}.

\smallskip

\noindent The parametric codimension is equal to
\begin{equation*}
\bar{M}(C,p)=\sum_{i=1}^{k_p} (m_i^p-1) +\omega_p-1=-L_p+\sum_{i=1}^{k_p} m_i^p \tag{7}.
\end{equation*}
For other equivalent formulae of $\bar{M}(C,p)$ see \cite{orza}.

\section{Global data}

\subsection{}\label{2.1} Let $\bp^N$, where  $ N=d(d+3)/2$,
 be the Hilbert scheme parametrising algebraic projective plane curves of degree $d$.
The locus $V_{d,g}$ in $\bp^N$ 
of reduced and irreducible  curves of degree $d$ and genus $g$ is irreducible 
(see e.g. \cite{Harris}). The locus of reduced and irreducible  curves of 
degree $d$ and genus $g$ having only nodes as singularities 
is smooth of dimension $3d-1+g$ and is dense in $V_{d,g}$. 

Recall that our goals is to  characterise the local
embedded topological types of local singular germs $(C,p_i)\subset (\bp^2,p_i)
$ which can be realized by
a projective curve $C$ of degree $d$.
Assume that a projective reduced plane curve $C$ of degree $d$
exists with fixed topological types
$S_1,\ldots,S_\nu$. Let $V(S_1,\ldots,S_\nu)$ denote the locally
closed  subscheme of reduced curves on $\bp^2$
having exactly $\nu$-singularities with topological types $S_1,\ldots,S_\nu$.

Greuel and Lossen in Theorem 3.6 of \cite{gl}
proved that the dimension of $V(S_1,\ldots,S_\nu)$ at $C$ satisfies
\begin{equation*}
\dim (V(S_1,\ldots,S_\nu),C)\geq C^2+1-p_a(C)-\tau^{es}(C)=\frac{d(d+3)}{2}-\tau^{es}(C), \tag{8}
\end{equation*}
where $\tau^{es}(C):=\sum_{i=1}^\nu \tau^{es}(C,p_i).$ The right hand side
\begin{equation*}
\text{expdim}(V(S_1,\ldots,S_\nu),C):=\frac{d(d+3)}{2}-\tau^{es}(C)
\end{equation*}
is called \emph{expected dimension} of $V(S_1,\ldots,S_\nu)$ at $C$.
The study of the locally closed subscheme $V(S_1,\ldots,S_\nu)$ of $\bp^N$
and its properties have been studied intensively in the literature (see for instance works by 
Artal Bartolo, Diaz, Greuel, Harris, Karras, Lossen, Shustin, Tannenbaum, Wahl, Zariski among many others).

\subsection{The action of $PGL(3,\bc)$}\label{2.2} Consider a reduced curve $C$ of
degree $d$ in $\bp^2$ and the action of the group $PGL(3):=
PGL(3,\bc)$ on the space $\bp^N$, which  parametrises plane curves of degree $d$. The orbit of a curve $C$ is a quasi-projective variety of dimension $\dim PGL(3)-\dim \text{Stab}_{PGL(3)}(C)$. For a general curve, the dimension of the orbit is $8$, that is, its stabiliser is $0$-dimensional.

Since the topological types $S_1,\ldots,S_\nu$ of singularities of $C$ and the
degree $d$ of $C$ remain constant under the $PGL(3)$-action, we consider the
\emph{virtual dimension} of $V(S_1,\ldots,S_\nu)$ at $C$ defined by
\begin{equation*}
\text{virtdim}(V(S_1,\ldots,S_\nu),C):=\text{expdim}(V(S_1,\ldots,S_\nu),C)-(8-\dim \text{Stab}_{PGL(3)}(C)).\tag{9}
\end{equation*}

Curves with small orbits have been studied and classified by P.~Allufi and
C.~Faber in \cite{AF}. According to this  classification,
 $C$ is always a configuration of rational  curves. Moreover,
 $C$ consists of irreducible components of the form below, with
arbitrary multiplicities.
We reproduce here their list together with
the dimension of the stabiliser $\text{Stab}_{PGL(3)}(C)$.
\begin{enumerate}
\item $C$ consists of a single line;
$\dim \text{Stab}_{PGL(3)}(C)=6$.
\item $C$ consists of 2 (distinct) lines; $\dim \text{Stab}_{PGL(3)}(C)=4$.
\item $C$ consists of 3 or more concurrent lines; $\dim \text{Stab}_{PGL(3)}(C)=3$.
\item $C$ is a triangle (consisting of 3 lines in general position);
$\dim \text{Stab}_{PGL(3)}(C)=2$.
\item $C$ consists of 3 or more concurrent lines, together with 1 other
(non-concur\-rent) line; $\dim \text{Stab}_{PGL(3)}(C)=1$.
\item $C$ consists of a single conic; $\dim \text{Stab}_{PGL(3)}(C)=3$.
\item $C$ consists of a conic and a tangent line; $\dim \text{Stab}_{PGL(3)}(C)=2$.
\item $C$ consists of a conic and 2 (distinct) tangent lines;
$\dim \text{Stab}_{PGL(3)}(C)=1$.
\item $C$ consists of a conic and a transversal line and may contain either
one of the tangent lines at the 2 points of intersection or both of them;
$\dim \text{Stab}_{PGL(3)}(C)=1$.
\item $C$ consists of 2 or more bitangent conics (conics in the pencil
$y^2+\lambda x z$) and may contain the line $y$ through the two points
of intersection as well as the lines $x$ and/or $z$, tangent lines to the
conics at the points of intersection; again, $\dim \text{Stab}_{PGL(3)}(C)=1$.
\item $C$ consists of 1 or more (irreducible) curves from the pencil
$y^b+\lambda z^a x^{b-a}$, with $b\ge 3$, and may contain the lines $x$
and/or $y$ and/or $z$; $\dim \text{Stab}_{PGL(3)}(C)=1$.
\item $C$ contains 2 or more
conics from a pencil through a conic and a
double tangent line; it may also contain that tangent line. In this case,
$\dim \text{Stab}_{PGL(3)}(C)=1$.
\end{enumerate}

\subsection{The Coolidge-Nagata conjecture.}

Let $C_1$ and $C_2$ be two curves in $\bp^2$. The pairs
$(\bp^2,C_1)$ and $(\bp^2,C_2)$ are \emph{birationally equivalent}
if there exist a birational map $\sigma:\bp^2\dashrightarrow
\bp^2$ such that the proper image of $C_1$ by $f$ coincides with
$C_2$. Traditionally,
the map $\sigma$ is called a \emph{Cremona transformation}.

\medskip

Let $\pi:V\to \bp^2$ be an embedded resolution of singularities of $C$,
let $\bar{C}$ be the strict transform of $C$ by $\pi$
and let $K_V$ be the canonical divisor of $V$.
One can show  that  $h^0(m(K_V+\bar{C}))$ and the Kodaira dimension
$\kappa(K_V+\bar{C},V)$ are birational
invariants of $C$ as plane curves (see \cite{Ii}, \cite{KM}).
Thus,  one defines
$\kappa[C]:=\kappa(K_V+\bar{C},V)$.

\medskip

Let $D$ be the reduced total preimage of $C$ by the embedded resolution $\pi$.
One can also show  that  $h^0(m(K_V+D))$ and the Kodaira dimension
$\kappa(K_V+D,V)$ are birational
invariants of the open surface $Y:=\bp^2\setminus C$ (see \cite{Ii}).
Thus one defines the logarithmic 
Kodaira dimension $\bar{\kappa}(Y):=\kappa(K_V+D,V)$ of $Y$.

\medskip

One of the  open problems regarding projective plane curves is the
following famous {\em Coolidge-Nagata problem/conjecture}, cf.
\cite{Coo} and \cite{Na}: {\em  every rational cuspidal curve can
be transformed by a Cremona transformation into a straight line}.

\medskip

In \cite{Coo}, J.I.~Coolidge proved:
\subsection{Theorem}
{\em A rational curve can be transformed into a straight line by a
Cremona transformation if and only if all the conditions for special adjoints
of any index are  incompatible.}

\vspace{2mm}

 A \emph{special adjoint of index $m$} is an effective
divisor in the complete linear system $|mK_V+\bar{C}|$. In fact
one can check that
 $\kappa[C]=-\infty$ if and only if $C$ has no special adjoints.
One has the following criterium due to
N.M.~Kumar and M.P.~Murthy (cf. Corollary 2.4 in \cite{KM}),
cf. also with Iitaka (Proposition 12 in  \cite{Ii}):

\subsection{Theorem}\label{KMI} {\em Let $C$ be an irreducible rational
plane curve. The following conditions are equivalent:
\begin{enumerate}
\item[a)] the curve $C$ can be transformed
into a straight line by a Cremona transformation,
\item[b)]  $\kappa[C]=-\infty$,
\item[c)] $|2K_V+\bar{C}|=\emptyset$,
\item[d)] $|2(K_V+\bar{C})|=\emptyset$.
\end{enumerate}}

\noindent
In \cite{KM} N.M.~Kumar and M.P.~Murthy also showed that a
sufficient condition is $\bar{C}^2\geq -3.$

\medskip

The Nagata-Coolidge problem  has been solved for cuspidal
rational plane curves with logarithmic Kodaira dimension of the complement
$\bar{\kappa}(\bp^2\setminus C)<2$ (using the classification
listed in \ref{classific}) and for
all known curves with $\bar{\kappa}=2$ (see below, cf. also with \ref{haha}).

\section{On rational curves}

From now on we are interested in rational projective plane curves.
Let $C$ be a reduced rational projective plane curve
 of degree $d$ in the complex projective plane  with singular points $\{p_j\}
_{j=1}^\nu$. Let $S_j$ be the topological type of the singularity $(C,p_j)$.

The main theme of this section is  Orevkov's conjecture \cite{Or}.
Since this conjecture is not true for a small number of curves
(with positive dimensional stabiliser ${Stab}_{PGL(3)}(C)$, all of
them with $\bar{\kappa}(Y)<2$), we correct the original version by
adding a contribution provided by this stabiliser; and  we also
present some equivalent  reformulations.

\subsection{\bf Virtual dimension conjecture.}\label{LuC}
\emph{
The virtual dimension {\em virtdim}$(V(S_1,\ldots,S_\nu),C)$
of reduced rational projective plane curve is  non-negative:}
\begin{equation*}
\text{expdim}(V(S_1,\ldots,S_\nu),C)-8+\dim \text{Stab}_{PGL(3)}(C))\geq 0.\tag{10}
\end{equation*}

\noindent This version was also conjectured independently by I.
Luengo.

\vspace{2mm}

Since $C$ is rational, if the multiplicity sequence of the singularity $(C,p)$ is
$[m_1^p,\ldots,m_k^p]$, then (the genus-formula reads as):
\begin{equation*}(d-1)(d-2)=\sum_{p\in\text{Sing}(C)} \sum_{i=1}^k m_i^p(m_i^p-1)=
\sum_{p\in\text{Sing}(C)} \left( \sum_{i=1}^k (m_i^p)^2-\sum_{i=1}^k m_i^p\right).\tag{11}
\end{equation*}
Eliminating  $d^2$ from  (9) and (5),  one gets
\begin{equation*}
\text{virtual dim}=3d-9-\sum_{p\in\text{Sing}(C)} \sum_{i=1}^k m_i^p +
\sum_{P\in\text{Sing}(C)} L_p + \dim \text{Stab}_{PGL(3)}(C).
\end{equation*}
Substituting (7) in this equality, one gets

\vspace{.5mm}

\begin{equation*}
\text{virtual dim}=3d-9-\sum_{p\in\text{Sing}(C)} \bar{M}(C,p)+
\dim \text{Stab}_{PGL(3)}(C). \tag{12}
\end{equation*}

In particular, $(10)$ is equivalent to

\vspace{.5mm}

\begin{equation*}
3d-9-\sum_{p\in\text{Sing}(C)} \bar{M}(C,p)+ \dim \text{Stab}_{PGL(3)}(C)\geq 0. \tag{13}
\end{equation*}

\subsection{\bf Orevkov's conjecture.}\label{OC}
S.~Orevkov in \cite{Or} conjectured the following inequality.

\medskip

\noindent {\em For a rational cuspidal curve}
\begin{equation*}
\sum_{p\in\text{Sing}(C)} \bar{M}(C,p) \leq 3d-9. \tag{14}
\end{equation*}

\noindent Orevkov's conjecture is stated for any rational curve (without
 restrictions). In such a case,
 the sum is over all irreducible branches at each singular point.
Note that (14) is not true for the curve $C$ defined by
 $x^2y+z^3=0$ since  $\bar{M}(C,p)=1$ at its  singular point.
Nevertheless,  $\dim \text{Stab}_{PGL(3)}(C)=1$ and (10) and (13)
hold. We believe that  the correct statement of the conjecture is
(10), or equivalently (13). Nevertheless, one can prove (cf. Lemma
\ref{AlF}) that for  irreducible, cuspidal, rational plane curve
with with $\bar{\kappa}(Y)=2$ the statements of  $(10)$ and $(14)$
are equivalent.

\subsection{\bf The $\cc^2$--Conjecture.}\label{CC}
Let $\pi:V\to
\bp^2$ be the minimal good embedded resolution of $C\subset
\bp^2$, and let   $\bar{C}$ be the strict transform of $C$ and
  $D=\pi^{-1}(C)$ be the reduced preimage of $C$
as above. One of the integers which plays a special
role in the classification problem is
the self-intersection of $\cc$ in $V$. It equals
\begin{equation*}
\cc^2=d^2-\sum_{p\in\text{Sing}(C)} \sum_{i=1}^k (m_i^p)^2 . \tag{15}
\end{equation*}
From $(10)$ and $(15)$ one gets
\begin{equation*}
3d=\cc^2+2+\sum_{p\in\text{Sing}(C)} \sum_{i=1}^k m_i^p.
\end{equation*}
Thus, via (12):
\begin{equation*}
\text{virtual dim}=\cc^2-7+ \sum_{p\in\text{Sing}(C)} L_p + \dim \text{Stab}_{PGL(3)}(C). \tag{16}
\end{equation*}
 Therefore,  inequality $(10)$ is equivalent to:
\begin{equation*}
\cc^2-7+ \sum_{p\in\text{Sing}(C)} L_p + \dim \text{Stab}_{PGL(3)}(C)\geq 0. \tag{17}
\end{equation*}

\section{Cuspidal rational curves and the rigidity conjecture.}\label{rigconj}

Let $C$ be an irreducible curve of degree $d$ in the complex projective plane.
One of the main invariants of such curves
is the logarithmic Kodaira dimension $\bar{\kappa}=\bar{\kappa}(Y)$,
where $Y:=\bp^2\setminus C$.
The following result of Wakabayashi  \cite{Wak} is crucial  in the
classification procedure.

\subsection{Theorem}\cite{Wak} \label{Wak} {\em Let $C$ be an
irreducible curve of degree $d$ in $\bp^2$.
\begin{enumerate}
\item If $g(C)\geq 1$ and $d\geq 4$ then $\bar{\kappa}(\bp^2\setminus C)=2$.
\item If $g(C)=0$ and $C$ has at least 3 cuspidal points then $\bar{\kappa}(\bp^2\setminus C)=2$.
\item If $g(C)=0$ and $C$ has at least 2 singular points and one of them is locally reducible then
$\bar{\kappa}(\bp^2\setminus C)=2$.
\item If $g(C)=0$ and $C$ has 2 cuspidal points then
$\bar{\kappa}(\bp^2\setminus C)\geq 0$.
\end{enumerate}}

\subsection{}\label{classific}
Recall that the  open surface
$Y:=\bp^2\setminus C$ is $\bq$-acyclic if and only if $C$ is a rational
cuspidal curve. From now on we assume that $C$ is a rational cuspidal plane
curve of degree $d$, in particular $C$ is irreducible.
Let  $(C,p_i)_{i=1}^\nu$ be the
collection of local plane curve  singularities, all of them locally
irreducible.

\vspace{2mm}

(a) If $\bar{\kappa}=-\infty$ then $\nu=1$ by \cite{Wak}. Moreover, all these
curves are classified by M.~Miyanishi and T.~Sugie \cite{MiS} (see also H. Kashiwara \cite{Kash}).
The family contains as an
important subfamily the Abhyankar-Moh-Suzuki curves (see~\cite{flmn}).

(b) The case $\bar{\kappa}=0$ cannot occur by a result of Sh.~Tsunoda
\cite{Tsu1}, see also Orevkov's article  \cite{Or}.

(c) If $\bar{\kappa}=1$ then by the above  result of Wakabayashi \cite{Wak}
one has $\nu\leq 2$. In the case  $\nu=1$, K. Tono provides  the
possible equations of the curves \cite{KeitaTono}. (Notice that
Tsunoda's  classification  in \cite{Tsu} is incomplete.). On the other hand, by
another  result of Tono \cite{To}, the case $\nu=2$ corresponds
exactly to the Lin-Zaidenberg bicuspidal rational plane curves.

\subsection{Lemma} \label{AlF} {\em If $C$ is an irreducible,
cuspidal, rational plane curve with with $\bar{\kappa}(Y)=2$ then
$\dim \text{Stab}_{PGL(3)}(C)=0$.}

\begin{proof}
According with the above classification of Aluffi and Faber, $C$
is a rational cuspidal plane curve with $\dim
\text{Stab}_{PGL(3)}(C)>0$ if and only if $C$ is a generic member
of the pencil $y^d+\lambda z^{a} x^{d-a}$ with $d\geq 3$ and
$(d,a)=1$. If $a=d-1$ then $C$ is of Abhyankar-Moh-Suzuki type
(with $\bar{\kappa}=-\infty$). Otherwise $a\ne d-1$ and $C$ is of
Lin-Zaidenberg type (with $\bar{\kappa}=1$).
\end{proof}

\subsection{} From a different point of view,
 one can classify triples $(d,a,b)$ such that there exists
a unicuspidal rational plane curve $C$ of degree $d$ whose singularity
has only one Puiseux pair of type
$(a,b)$,  where $1<a<b$. Let $\{\varphi_j\}_{j\geq 0}$ denote
the Fibonacci numbers $\vp_0=0$, $\vp_1=1$,
$\vp_{j+2}=\vp_{j+1}+\vp_j$.

\subsection{Theorem} \cite{class} {\em The Puiseux pair $(a,b)$ can be
realized by a unicuspidal rational curve of degree $d$ if and only
if $(d,a,b)$ appears in the following list.

 \vspace{1mm}

(a) $(a,b)=(d-1,d)$;

(b) $(a,b)=(d/2, 2d-1)$, where $d$ is even;

(c) $(a,b)=(\vp_{j-2}^2,\vp_j^2)$ and $d=\vp_{j-1}^2+1=\vp_{j-2} \vp_{j}$, where $j$ is odd
and $\geq 5$;

(d) $(a,b)=(\vp_{j-2},\vp_{j+2})$ and $d=\vp_{j}$, where $j$ is odd and
$\geq 5$;

(e) $(a,b)=(\vp_4,\vp_8+1)=(3,22)$ and $d=\vp_6=8$;

(f) $(a,b)=(2\vp_4,2\vp_8+1)=(6,43)$ and $d=2\vp_6=16$.}

\vspace{2mm}

In the first four cases $\bar{\kappa}=-\infty$ and they can be
realized by some particular curves which appear  in Kashiwara's
classification \cite{Kash}. The last two sporadic cases have
$\bar{\kappa}=2$ and were found by Orevkov and Artal-Bartolo, cf.
\cite{Or}.

\subsection{} Another aproach is to classify rational cuspidal  curves $C$ such that the highest  multiplicity of the singular
points $m$ is close to the degree $d$. Flenner and Zaidenberg classified the curves with $m=d-2$  in \cite{FZ1} and  $m=d-3$ in \cite{FZ2}. The case 
$m=d-4$ is partially solved by Fenske~\cite{Fenske1}. Note that $m$ can not be too small
because in \cite{MS} it is proved that $d < 3m$ solving a conjecture of Yoshihara~\cite{Yos}. Let $\alpha = {(3 + \sqrt{5}) \over2}$, Orevkov~\cite{Or} gives two families of curves with  ${\alpha}m < d$ and conjectured that those families gives the only curves verifyng ${\alpha}m < d.$

\subsection{The rigidity conjecture of Flenner and Zaidenberg} \
Let $Y$ be an $\bq$-acyclic affine surface, and fix one of its
 `minimal logarithmic compactifications' $(V,D)$. This means that $V$
is a smooth projective surface with a normal crossing divisor $D$, such that
$Y=V\setminus D$, and $(V,D)$ is minimal with these properties.

The sheaf of the logarithmic tangent vectors $\Theta_V\langle
D\rangle$ controls the  deformation theory of the pair $(V,D)$,
cf. \cite{FZ0}. E.g., $H^0(V,\Theta_V\langle D\rangle )$ is the
set of infinitesimal automorphisms, $H^1(V,\Theta_V\langle
D\rangle )$ is the space of infinitesimal deformations, and
$H^2(V,\Theta_V\langle D\rangle )$ is the space of obstructions.
Iitaka showed in \cite{Ii2} that if $\bar{\kappa}(Y)=2$ then the
automorphism group of the surface $Y$ is finite (this also
provides a different proof of \ref{AlF}). Therefore
$h^0(\Theta_V\langle D\rangle )=0$. In \cite{FZ0,FZ1,Z} Flenner and
Zaidenberg  conjectured the following

\medskip

\noindent {\bf  Rigidity conjecture:}
 {\em Every $\bq$-acyclic affine surfaces $Y$
with logarithmic Kodaira dimension $\bar{\kappa}(Y)=2$
is rigid and has unobstructed deformations. That is,
\begin{equation*}
h^1(\Theta_V\langle D\rangle )=0 \quad \quad \text{and} \quad
\quad h^2(\Theta_V\langle D\rangle )=0. \tag{18}
\end{equation*}
In particular,  the Euler characteristic $\chi(\Theta_V\langle D\rangle )=
h^2(\Theta_V\langle D\rangle )-h^1(\Theta_V\langle D\rangle )$ must vanish.
}

\medskip

In~\cite{FZ0},~\cite{FZ1} and~\cite{FZ2} the conjecture was verified for most of the known examples. In~\cite{FZ0} unobstructedness
was proved for all $\bq$-acyclic surfaces of non log-general type. In~\cite{ZO} it is proved that a  rigid rational cuspidal curve has at most 9 cusps.

This can be applied in our situation as follows. Consider a
projective curve $C$, and write $Y:=\bp^2\setminus C$.
The $\bq$-acyclicity of $Y$ is equivalent to the fact that
$C$ is rational and cuspidal. For $V$ one can take the minimal embedded
resolution of the pair $(\bp^2,C)$.

The conjecture for $Y=\bp^2\setminus C$  implies
the projective rigidity of  the curve $C$. This means that
every equisingular deformation of $C$ in $\bp^2$ would be
projectively equivalent to $C$. Thus
$V(S_1,\ldots,S_\nu)$ has expected dimension 8 (see Section 4).

\medskip

In Corollary 2.5 of \cite{FZ1}, Flenner and Zaidenberg  show that
for any cuspidal rational plane curve
\begin{equation*}
\chi(\Theta_V\langle D\rangle )=K_V(K_V+D)=-3(d-3)+\sum_{p\in \text{Sing}(C)}
\bar{M}(C,p). \tag{19}
\end{equation*}
By $(12)$ and Lemma \ref{AlF} then
\begin{equation*}
\text{virtual dim}=-\chi(\Theta_V\langle D\rangle ). \tag{20}
\end{equation*}

\medskip

\noindent
The vanishing of $\chi(\Theta_V\langle D\rangle )$ implies any of the
equivalent equalities $(10)$,  $(13)$ or  $(17).$
On the other hand, if $(10)$,  $(13)$ or  $(17)$ hold,  then
\begin{equation*}
\chi(\Theta_V\langle D\rangle )\leq 0. \tag{21}
\end{equation*}

\medskip

\subsection{Proposition}\label{key}  [Tono]  \ {\em For cuspidal
rational plane curves with $\bar{\kappa}=2$ the following
inequality holds}
\begin{equation*}
\chi(\Theta_V\langle D\rangle )\geq 0. \tag{22}
\end{equation*}

\medskip

\begin{proof}
(22) follows from the article  \cite{Tono2} of  K.~Tono in the following way.
F.~Sakai in \cite{Sakai} introduce the invariant $\gamma_2:=h^0(2K_V+D)$.
Lemma 4.1 in \cite{Tono2} states that if
the pair $(V,D)$ satisfies the following three conditions (for details see
[loc. cit.])
\begin{enumerate}
\item[(A1)] $\bar{\kappa}(V\setminus D)=2$,
\item[(A2)] $(V,D)$ is almost minimal, and
\item[(A3)] $D$ contains neither a rod consisting of $(-2)$-curves nor a fork consisting of $(-2)$-curves,
\end{enumerate}
then $$\gamma_2=K_V(K_V+D)+\frac{D(D+K_V)}{2}+\chi(\mathcal{O}_V).$$
(The main point here is that by a vanishing theorem $h^1(2K_V+D)=0$, by an
 easy argument  $h^2(2K_V+D)=0$ too, 
hence $\gamma_2=\chi(2K_V+D)$ can be computed by Riemann-Roch.)

One can check that in our case the minimal embedded resolution
satisfies these  conditions.  Moreover, $\chi(\mathcal{O}_V)=1$ and
 (since $D$ is a rational tree, the adjunction formula implies)
$K_VD+D^2=-2$. Thus   $\gamma_2=K_V(K_V+D)$.
Therefore, via (19), one has:
\begin{equation*}
\chi(\Theta_V\langle D\rangle )= h^0(2K_V+D)\geq 0.
\end{equation*}
\end{proof}

\subsection{Corollary}\label{haha} {\em Let $C$ be an irreducible,
cuspidal, rational plane curve with $\bar{\kappa}(\bp^2-C)=2$. The
following conditions are equivalent:
\begin{enumerate}
\item[(i)] $\chi(\Theta_V\langle D\rangle )=0$,
\item[(ii)] $\text{virtdim}(V(S_1,\ldots,S_\nu),C)\geq 0$, i.e. $(10)$ holds, where
$S_j$ is the topological type of the corresponding uni-branch
singularity $(C,p_j)$.
\item[(iii)] $\chi(\Theta_V\langle D\rangle )\leq 0$.
\end{enumerate}
In such a case, the curve $C$ can be  transformed by a Cremona
transformation of $\bp^2$ into a straight line (i.e.,
the Coolidge-Nagata problem has a positive answer).}

\begin{proof}
$(i)\Rightarrow (ii)\Rightarrow (iii)$ follows from $(20)$ and
$(21)$. $(iii)\Rightarrow (i)$ follows from \ref{key} and $(21)$.
Finally, the characterisation  \ref{KMI} shows that $C$ can be
transform into a straight line by a Cremona  transformation.
Indeed, $h^0(2K_V+D)=\chi(\Theta_V\langle D\rangle) =0$, but
${\mathcal O}_V(2K_V+\bar{C}) $ is a subsheaf of ${\mathcal
O}_V(2K_V+D)$, hence  $h^0(2K_V+\bar{C})=0$ as well.
\end{proof}

\section{The semigroup distribution  property}

\subsection{}\label{conjecture1} The characterisation  problem of the
realization of prescribed topological types of singularities has a
long and rich history providing many interesting compatibility
properties connecting local invariants of the germs
$\{(C,p_i)\}_i$ with some global invariants of $C$ ---\, like its
degree, or the log-Kodaira dimension of $\bp^2\setminus C$, etc.
(For a --- non-complete --- list of some of these restrictions, see
e.g. \cite{flmn,class}.)

In \cite{flmn} we proposed a new compatibility property ---\,
valid for rational cuspidal curves $C$. Its formulation is
surprisingly very elementary. Consider a collection
$(C,p_i)_{i=1}^\nu$ of locally irreducible plane curve
singularities (i.e. cusps), let $\Delta_i(t)$ be the
characteristic polynomial of the monodromy action associated with
$(C,p_i)$, and $\Delta(t):=\prod_{i}\Delta_i(t)$. Its degree is
$2\delta$, where $\delta$ is the sum of the delta-invariants $\delta(C,p_i)$ of
the singular points. Then $\Delta(t)$ can be written as
$1+(t-1)\delta+(t-1)^2Q(t)$ for some polynomial $Q(t)$. Let $c_l$
be the coefficient of $t^{(d-3-l)d}$ in $Q(t)$ for any
$l=0,\ldots, d-3$.

\vspace{2mm}

\subsection{Conjecture A}\label{CA}  \cite{flmn} {\em Let
$(C,p_i)_{i=1}^\nu$ be a
 collection of local plane curve  singularities, all of them locally
irreducible,  such that
$2\delta=(d-1)(d-2)$ for some integer $d$. If $(C,p_i)_{i=1}^\nu$
can be realized as the local singularities of a degree $d$
(automatically rational and cuspidal)
projective plane curve then}
\begin{equation*}c_l\leq (l+1)(l+2)/2 \ \ \mbox{{\em for all}
 \ $l=0,\ldots, d-3$.}
\end{equation*}

In fact, the integers $n_l:=c_l-(l+1)(l+2)/2$ \
are symmetric: $n_l=n_{d-3-l}$; and $n_0=n_{d-3}=0$.
We also mention that examples with strict inequality occur,
cf. \cite{flmn} (in all these examples known by the authors
$\nu>1$).

The main result of \cite{flmn} is :

\subsection{Theorem}\cite{flmn}  If {\em $\bar{\kappa}(\bp^2\setminus C)$ is $\leq 1$, then
the above conjecture A is true (in fact with $n_l=0$).}
 \vspace{2mm}

\noindent There is an additional  surprising phenomenon in the
above conjecture. Namely, in the {\em unicuspidal}  case one can
show the following.

\subsection{Proposition}\label{nuone}\cite{flmn}  {\em
If $\nu=1$ then \, $c_l\geq (l+1)(l+2)/2$ \, for \, $0\leq l\leq d-3$.}

\vspace{2mm}

Therefore, conjecture \ref{CA} in this case can be reformulated as
follows:

\subsection{Conjecture B1}\label{CB1} {\em With the notations of
\ref{CA}, if $\nu=1$, then $n_l=0$ for all $l=0,\ldots, d-3$, that is}
\begin{equation*}c_l=(l+1)(l+2)/2 \ \ \mbox{{\em for all}
 \ $l=0,\ldots, d-3$.}\end{equation*}

In fact, if $\nu=1$, we can do more.
 Recall that the characteristic
polynomial $\Delta$ of $(C,p)\subset (\bp^2,p)$ is a complete
(embedded) topological invariant of this germ, similarly as the
semigroup $\Gamma_{(C,p)}\subset \bn$. In the next discussion we
will replace $\Delta $ by $\Gamma_{(C,p)}$. Recall that the
semigroup $\Gamma_{(C,p)}\subset \bn$ consists of all possible
intersection multiplicities $I_p(C,h)$ at the point $p$ for all
$h\in {\mathcal O}_{(\bc^2,p)}$.

Hence, one can reformulate  conjecture B1 in terms of the
semigroup of the germ $(C,p)$ and the degree $d$. It turns
out that the of vanishing of the coefficients
$n_l$ (as in B1.) is replaced  by a very precise and mysterious distribution
of the elements of the semigroup with
respect to the intervals $I_l:=(\,(l-1)d,ld\,]$:

\subsection{Conjecture B2}\label{CB2} 
 {\em Assume that $\nu=1$. Then for   any $l>0$,  the interval $I_l$ 
contains exactly $\min\{l+1,d\}$  elements from the  semigroup $\Gamma_{(C,p)}$.}

\vspace{2mm}

In other words,  for every rational unicuspidal plane curve $C$ of
degree $d$, the above conjecture is equivalent to the
identity
\begin{equation*}
D(t)\equiv 0,\tag{$DP$}
\end{equation*}
where:
\begin{equation*}
D(t):= \sum_{k \in \Gamma_{(C,p)}} t^{\lceil k/d\rceil}
-\Big(1+2t+\cdots+
(d-1)t^{d-2}+d(t^{d-1}+t^d+t^{d+1}+\cdots)\Big).
\end{equation*}

For the equivalences of conjectures B1 and
B2, see \ref{top}. Here we only mention a key relation between
the coefficients $c_l$ and the semigroup $\Gamma_{(C,p)}$.

 First,
consider the identity (cf. \cite{CDG}) $\Delta(t)=(1-t)\cdot
L(t)$, where $L(t) =\sum_{k\in \Gamma_{(C,p)}} t^k$ is the
Poincar\'e series of $\Gamma_{(C,p)}$. Write
$\Delta(t)=1-P(t)(1-t)$ for some polynomial $P(t)$, then
$L(t)+P(t)=1/(1-t)=\sum_{k\geq 0}t^k$. In particular,
$P(t)=\sum_{k \in \bn\setminus \Gamma_{(C,p)}} t^k$. Then
$$Q(t)=\frac{P(t)-\delta}{t-1}=\sum_{k\not\in \Gamma_{(C,p)}}
\frac{t^k-1}{t-1}=
\sum_{k\not\in \Gamma_{(C,p)}}
(1+t+\cdots +t^{k-1}).$$
Hence $c_l=\#\{
k\not\in \Gamma_{(C,p)} \, :\, k>(d-3-l)d\}$.
Since $k\in \Gamma_{(C,0)}$ if and only if
$\mu-1-k \not\in \Gamma_{(C,p)}$  for any $0\leq k\leq \mu-1$, one gets
$c_l=\#\{k\in \Gamma_{(C,p)}\,;\, k\leq ld\}$.

\subsection{}
The following equivalent formulation was suggested by A. Campillo.

\subsection{Theorem} {\em Let $C$ be a unicuspidal rational plane curve of
degree $d$.
The curve $C$ satisfies the semigroup compatibility property
$(DP)$ (i.e. conjectures B1 and/or B2) if and only if the elements
of the semigroup $\Gamma_{(C,p)}$ in $[0,ld]$ are realized by
projective (possibly non-reduced) curves of degree $l$ for $l\leq
d-3$.}

\begin{proof}
The proof of the `if' part is easy.
For the `only if' part fix a projective coordinate system  $[X:Y:Z]$
such that the affine chart $Z\not=0$ contains  the singular point $p$.
Let $V$ be the vector space of polynomials of degree $l$
in variables $(X/Z,Y/Z)$. Its dimension  is $N:=(l+1)(l+2)/2$,
which, in fact, equals the number of elements of the semigroup in the interval
$[0,ld].$ Denote these elements by $0=s_1,...,s_{N}$, ordered in an increasing way.

Consider the decreasing filtration of vector spaces $V_1\supset V_2\supset \cdots \supset  V_N,$
defined by
$$
V_i:=\{f\in V: I_p(C,f)\geq s_i\}.
$$
First, we verify that $\dim (V_i/V_{i+1})$ is at most $1$.
Indeed, assume that $I_p(C,f_i)=I_p(C,f_2)=I$. Let $n:(\bc,0)\to (C,p)$
be the normalisation of $(C,p)$, and write $f_i\circ n(t)=
a_it^I+\cdots$ with $a_i\not=0$, for $i=1$ and 2. Then $I_p(C,a_2f_1-
a_1f_2)>I$. Since there is no semigroup element between  $s_i$ and
$s_{i+1}$, the inequality $\dim (V_i/V_{i+1})\leq 1$ follows.

Next, notice that to prove the theorem  it is enough to show that each
dimension $\dim(V_i/V_{i+1})$ is exactly $1$.

But, if
$\dim(V_i/V_{i+1})=0$ for some $i$ then $\dim(V_N)$ is at
least $2$.
Since for any $f\in V_N$ one has $I_p(C,f)\geq s_N$,  $\dim(V_N)\geq 2$
would imply (by similar argument as above) the existence of an
$f\in V_N$ with $I_p(C,f)>s_N$. Since $I_p(C,f)$ is  an element of
the semigroup and the last element of the semigroup in the interval $[0,ld]$ is
$s_N$, we get that $I_p(C,f)>ld$,
which contradicts the irreducibility of $C$ by
B\'ezout Theorem.
\end{proof}

\subsection{A counterexample to an `extended' version}

In \cite{flmn} we formulated  the following conjecture,
as an extension of the conjecture B2. to an `if and only if'
statement.

\subsection{`Conjecture' C} {\em The local topological type $(C,p)\subset
(\bp^2,p)$ can be realized by a degree $d$ unicuspidal rational
curve if and only if the property $(DP)$ is valid. }

\medskip

In the sequel we present a counterexample to the `if'
part (i.e. to the `extension').

\vspace{2mm}

If the germ $(C,0)$ has $g$  Newton pairs
$\{(p_k,q_k)\}_{k=1}^g$
with gcd$(p_k,q_k)=1$, $p_k\geq 2$ and $q_k\geq 1$
(and by convention, $q_1>p_1$),
define the integers $\{a_k\}_{k=1}^g$ by
$a_1:=q_1$ and $a_{k+1}:=q_{k+1}+p_{k+1}p_ka_k$ for  $k\geq 1$.
Then its  Eisenbud-Neumann splice diagram decorated by
the numerical data $\{(p_k,a_k)\}_{k=1}^g$ has the following shape \cite{en}:

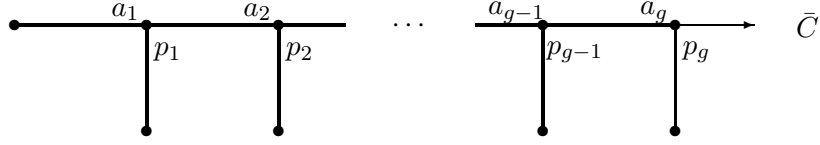
\begin{picture}(400,65)(0,10)
\put(50,60){\circle*{4}}
\put(100,60){\circle*{4}}
\put(150,60){\circle*{4}}
\put(250,60){\circle*{4}}
\put(300,60){\circle*{4}}
\put(100,20){\circle*{4}}
\put(150,20){\circle*{4}}
\put(250,20){\circle*{4}}
\put(300,20){\circle*{4}}
\put(350,60){\makebox(0,0){$\bar{C}$}}
\put(92,65){\makebox(0,0){$a_1$}}
\put(142,65){\makebox(0,0){$a_2$}}
\put(240,65){\makebox(0,0){$a_{g-1}$}}
\put(292,65){\makebox(0,0){$a_g$}}
\put(108,50){\makebox(0,0){$p_1$}}
\put(158,50){\makebox(0,0){$p_2$}}
\put(262,50){\makebox(0,0){$p_{g-1}$}}
\put(308,50){\makebox(0,0){$p_g$}}
\put(200,60){\makebox(0,0){$\cdots$}}
\put(50,60){\framebox(125,0){}}
\put(225,60){\framebox(75,0){}}
\put(100,20){\framebox(0,40){}}
\put(150,20){\framebox(0,40){}}
\put(250,20){\framebox(0,40){}}
\put(300,20){\framebox(0,40){}}
\put(300,60){\vector(1,0){30}}
\end{picture}

Consider now the local singularity  whose Eisenbud-Neumann splice
diagram is decorated by two pairs $(p_1,a_1)=(2,7)$ and
$(p_2,a_2)=(4,73)$.
A local equation for such singularity can be
$(x^2-y^7)^4+x^{33}y=0$. Its multiplicity sequence is
$[8_3,4_6,1_4]$. A minimal set of generators of its semigroup
$\Gamma_{(C,p)}$ is given by $\langle 8,28,73\rangle $. Its Milnor
number is $16\cdot 15$, hence a possible unicuspidal plane curve
$C$ of degree $17$ might
 exist with such local singularity.
Moreover the distribution  property $(DP)$ of the semigroup is
also satisfied. Nevertheless, such a  curve $C$ does not exist. To
prove this, one can either use Cremona transformations to
transform $C$ into another
 curve for which one sees that it does not exist, or one uses
Varchenko's semi-continuity criterium for the spectrum of the
singularity \cite{Var1,Var2}. Here we will follow the second
argument.

The spectrum of the irreducible singularity $(C,0)$ can be computed
from the Newton pairs of the singularity.
The forth author provided such a formula in \cite{Sp}. It is convenient to 
consider the spectrum $Sp(C,0)=\sum_r\, n_r(r)$ as an element of
$\Z[\bq\cap (0,2)]$. We write $Sp_{(0,1)}(C,0)$ for the collection of 
spectral elements situated in the interval $(0,1)$. 

\subsubsection{}\label{Sp}{\bf Theorem.} {\em If the irreducible germ $(C,0)$ has
$g$ Newton pairs $\{(p_k,q_k)\}_{k=1}^g$ then
$$
Sp_{(0,1)}(C,0)=\sum_{k=1}^g S_k \quad \text{where} \quad S_k=\sum
\left(\frac{i/a_k+j/p_k+t}{p_{k+1}p_{k+2}\cdots p_g}\right),
$$
where the second sum is over $0<i<a_k,\ 0<j<p_k,\ i/a_k+j/p_k<1$
and $0\leq t\leq p_{k+1}p_{k+2}\cdots p_g-1$ (if $k=g$ then
$S_g=\sum (l/a_g+k/p_g)$ where the sum is over $0<l<a_g, \
0<k<p_g,\ l/a_g+k/p_g<1$).}

\vspace{2mm}

If the local singular type $\{(C,p)\}$ can be realized by a degree
$d$ plane curve $C$, then $(C,p)$ is in the deformation of the
`universal' plane germ $(U,0):=(x^d+y^d,0)$. In particular, the
collection of all spectral numbers $Sp(C,p)$ of the local plane
curve singularity $(C,p)$ satisfies the semi-continuity property
compared with the spectral numbers of $(U,0)$ for any interval
$(\alpha,\alpha+1)$. Since the spectral numbers of $(U,0)$ are of
type $l/d$, the semi-continuity property for intervals
$(-1+l/d,l/d)$ ($l=2,3,\ldots , d-1$)  reads as follows:
\begin{equation*}
\#\{\alpha\in Sp(C,p)\ :\ \alpha<l/d\}\leq (l-2)(l-1)/2. \tag{23}
\end{equation*}
In our case, for  $d=17$ and  $l=12$, using Theorem \ref{Sp} we get
\begin{equation*}
\#\{\alpha\in Sp(C,p)\ :\ \alpha<12/17\}-(12-2)(12-1)/2=1,
\end{equation*}
which contradicts (23).
Thus the rational unicuspidal plane curve $C$ of degree $17$
with such singularity cannot exist.

\vspace{1mm}

Thus, in the realization problem, the above case
$(p_1,a_1;p_2,a_2;d)$ cannot be eliminated by the semigroup
distribution property ($DP$), but it can be eliminated by the semi-continuity of
the spectrum. However it is not true that the semi-continuity implies
($DP$). For a more precise discussion see \cite{flmn}.

\section{The semigroup compatibility property and surface
singularities}

\subsection{Superisolated singularities}\label{egy}
 The theory of normal surface singularities
(in fact, of isolated hypersurface surface singularities)
`contains' in a canonical way the theory of complex projective
plane curves via the family of  {\em superisolated} singularities.
These singularities were introduced by  the second author in
 \cite{Ignacio}, see also \cite{alm} for a survey on them.  
A hypersurface singularity $f:(\bc^3,0)\to (\bc,0)$, $f=f_d+l^{d+1}$ (where $f_d$
is homogeneous of degree $d$ and $l$ is linear) is superisolated
if the projective plane curve $C:=\{f_d=0\}\subset \bp^2$ is
reduced, and none of its singularities $\{p_i\}_{i=1}^{\nu}$ is
situated on  $\{l=0\}$. The equisingular type of $f$ depends only
on $f_d$, i.e. only on the projective curve $C\subset \bp^2$. In
particular, all the invariants (of the equisingular type) of $f$
can be determined from the invariants of the pair $(\bp^2,C)$.

In the next discussion we follow \cite{flmn,[59]}. There is a
standard procedure which provides the plumbing graph of the link
$M$ of $f$ from the embedded resolution graphs of $(C,p_i)$'s and
the integer $d$. The point is that the link $M$ is a rational
homology sphere if and only if  $C$ is rational and cuspidal. In
this section, we will assume that these conditions are satisfied.
Let $\mu_i=\mu(C,p_i)$ and $\Delta_i$ be the Milnor number and the
characteristic polynomial of the local plane curve singularities
$(C,p_i)$. Set $2\delta:=\sum_i\mu_i$, $\Delta:=\prod_i
\Delta_i$, and $\bar{\Delta}(t):=t^{-\delta}\Delta(t)$.

Let $(V,D)$ be the minimal embedded resolution of the pair
$(\bp^2,C)$ as above. The minimal plumbing graph of $M$ (or,
equivalently, the minimal good resolution graph of the surface singularity 
$\{f=0\}$) can
be obtained from the dual graph of $D$ by decreasing   the
decoration (self-intersection) of $\bar{C}$ by $d(d+1)$.   In the
language of topologists, if $C$ is unicuspidal ($\nu=1$), then
$M=S_{-d}^3(K)$ (i.e. $M$ is obtained via surgery of the 3-sphere
$S^3$ along $K$ with surgery coefficient $-d$), where $K\subset
S^3$ is the local knot of $(C,p)$. One can also verify that
$H_1(M,\Z)=\Z_d$.

Another topological invariant of $f$ is the following one. Let $Z\to
(\{f=0\},0)$ be the minimal good resolution, $K_Z$ be the
canonical divisor of $Z$ and $\#$ the number of irreducible components of the
exceptional divisor (which equals the number of irreducible
components of $D$). Then $K_Z^2+\#$ is a well-defined invariant of
$f$, which, in fact, can be computed from the link  $M$ (or, from
its graph) as well. In our case, surprisingly, in this invariant
of the link $M$ all the information about the local types
$(C,p_i)$ are lost: $K_Z^2+\#=1-d(d-2)^2$, it depends only on $d$.

The same is true for the Euler characteristic $\chi(F)$, or for the
signature $\sigma(F)$ of the Milnor fiber $F$ of $f$, or about the
geometric genus $p_g$ of $f$.
 In fact, it is well-known that  for any hypersurface  singularity,
 any of $p_g$, $\sigma(F)$ and
$\chi(F)$ determines the remaining two modulo $K_Z^2+\#$. E.g.,
one has the relation:
\begin{equation*}
8p_g+\sigma(F)+K_Z^2+\#=0. \tag{24}\end{equation*} In our case,
for the superisolated singularity $f$, one has
$p_g=d(d-1)(d-2)/6$, hence the smoothing  invariants $\chi(F)$ and
$\sigma(F)$ depend only on the degree $d$.

\subsection{} For a normal surface singularity with rational homology sphere link
(and with some additional analytic restriction, e.g. complete
intersection or Gorenstein property)  there is a subtle connection
between the Seiberg-Witten invariants of its link $M$ and some
analytic/smoothing invariants. The hope is that the geometric
genus (or, equivalently, $\chi(F)$ or $\sigma(F)$, see (24) and
the discussion nearby), can be recovered from the link. The
starting point is an earlier conjecture of Neumann and Wahl
\cite{NW}:

 \subsubsection{} {\em For any isolated complete intersection whose link
$M$ is an integral homology sphere we have the equality $\sigma(F)=8\lambda(M)$, where $\lambda(M)$ is the Casson invariant 
of the link.}

\vspace{1mm}

Notice that the link of a hypersurface superisolated singularity
is never an integral homology sphere. The generalised conjecture,
applied to rational homology spheres  (Conjecture $SWC$ below) was
proposed by the forth author in a joint work with L.  Nicolaescu
in \cite{[51]} involving the Seiberg-Witten invariant of the link.
It was verified for rather large number of non-trivial special
families (rational and elliptic singularities, suspension
hypersurface singularities $f(x,y)+z^n$ with $f$ irreducible,
singularities with good $\bc^*$ action)
\cite{[51],[52],[55],[61],[69]}.  But the last three authors of
the present article have shown in \cite{[59]} that the conjecture
fails in general. The
counterexamples were provided exactly by superisolated
singularities and/or their universal abelian covers, see also Stevens paper
\cite{Ste} where he computes explicit equations for the universal abelian covers.
Nevertheless,
in the next paragraph 
we will recall this conjecture (in its original form), since this
have guided us to the semigroup compatibility property, and we
believe that it hides a deep mathematical substance (even if at
this moment it is not clear for what family we should expect its
validity).

Let $\ssw_M(can)$ be the Seiberg-Witten invariant
of the link $M$ associated with the canonical $spin^c$ structure
(this is induced by the complex structure of $\{f=0\}\setminus
\{0\}$, and it can be identified combinatorially from the graph of
$M$; in this article we will not discuss the invariants associated
with the other $spin^c$ structures).

 \subsection{`Conjecture' SWC}\label{SWC}
  \cite{[51]} {\em For a $\bq$-Gorenstein surface singularity
whose link $M$ is a rational homology sphere one has
$$\ssw_M(can)-(K_Z^2+\#)/8 =p_g.$$
In particular, if the singularity is Gorenstein and admits a
smoothing, then $-\ssw_M(can)=\sigma(F)/8$ (cf. (24)). }

\vspace{2mm}

 If $M$ is an integral  homology sphere then
$\ssw_M(can)=-\lambda(M)$. If $M$ is a rational homology sphere then by
a result of Nicolaescu \cite{Nico5},
$\ssw_M(can)=\et_M-\lambda(M)/|H_1(M,\Z)|$,  where $\lambda(M)$ is the
Casson-Walker invariant of $M$ (normalised as in \cite{Lescop}),
 and $\et_M$ denotes the
 sign refined Reidemeister-Turaev torsion  (associated with the
 canonical $spin^c$ structure) \cite{Tu5}.

\subsection{} In our present situation, when $M$ is the link of a
superisolated singularity $f$, one shows, cf. \cite{[59]}  (using
the notations of \ref{egy}), that
\begin{equation}
\et_M=\frac{1}{d}\sum_{\xi^d=1\not=\xi}\
\frac{\Delta(\xi)}{(\xi-1)^2} \ \ \ \ \ \mbox{and } \ \ \ \ \
\lambda(M)= -\frac{ \bar{\Delta}(t)''(1)}{2} +
\frac{(d-1)(d-2)}{24}. \tag{25}\end{equation} Therefore, since
$p_g$ and $K_Z^2+\#$ depend only on $d$, the $SWC$ imposes serious
restriction on the local invariant $\Delta$. This condition, for
some cases when the number of singular points of $C$ is $\geq 2$,
is not satisfied (hence $SWC$ fails, cf. \cite{[59]});
nevertheless, as we will see, the $SWC$ identity in the
unicuspidal case is equivalent with Conjecture B2 of section 6
about the distribution property of the semigroup. In order to
explain this, let us  {\em assume that $C$ is unicuspidal}, and consider
(motivated by (25))
\begin{equation*}R(t):= \frac{1}{d}\sum_{\xi^d=1}
\frac{\Delta(\xi t)}{(1-\xi t)^2}-\frac{1-t^{d^2}}{(1-t^d)^3}.
\end{equation*}
Similarly,
\begin{equation*}
N(t):=\sum_{l=0}^{d-3}\Big(
c_l-\frac{(l+1)(l+2)}{2}\Big)t^{d-3-l};\ \ \ \mbox{and}\ \ D(t):=
\sum_{k\in\Gamma_{(C,p)}} t^{\lceil k/d\rceil}-\frac{1-t^d}{(1-t)^2}.
\end{equation*}
Notice that this $D(t)$ agrees with the one defined in ($DP$),
section 6. In \cite{flmn} the following facts are verified:
\begin{equation*}
R(t)=D(t^d)/(1-t^d)=N(t^d).\tag{26}\end{equation*}

\begin{equation*}
\mbox{$N(t)$ (hence $R(t)$ too)  has non-negative
coefficients.}\tag{27}\end{equation*} \begin{equation*}
R(1)=\ssw_M(can)-\frac{K^2+\#}{8}-p_g.\tag{28}\end{equation*}

Therefore, in this case, we have the equivalence of the `Seiberg-Witten
invariant conjecture' with the `semigroup distribution property':

\subsection{Theorem}\label{top} {\em Assume that $C$ is unicuspidal and rational (that is, $\nu=1$). Then the following facts are
equivalent:

(a) $R(1)=0$, i.e. Conjecture $SWC$ (\ref{SWC}) is true (for the above germ $f$);

(b) $R(t)\equiv 0$;

(c) $N(t)\equiv 0$, i.e. Conjecture B1 (\ref{CB1}) is true;

(d) $D(t)\equiv 0$, i.e. Conjecture B2 (\ref{CB2}) is true.}

\section{The semigroup distribution property and Heegaard Floer homology}

\subsection{} The presentation of this section is based on
some recent results of the forth author in \cite{[61],[63],[65]}.
In the sequel we assume that $C$ is {\em unicupidal},
and we keep the notations of the previous section.

There is another way  to compute the Seiberg-Witten invariant of
the link $M$ via its Heegaard-Floer homology.   For any oriented
rational homology 3-sphere $M$ the Heegaard Floer homology
$HF^+(M)$ was introduced by Ozsv\'ath and Szab\'o in \cite{OSz}
(cf. also with their long list of articles). $HF^+(M)$ is a
$\Z[U]$-module with compatible $\Q$-grading. Moreover, $HF^+(M)$
has a natural direct sum decomposition (compatible with the
$\Q$-grading) corresponding to the $spin^c$-structures of $M$: In
this article we write $HF^+(M,can)$ for the Heegaard-Floer
homology associated with the canonical $spin^c$ structure.

For some (negative definite) plumbed rational homology 3-spheres
$M$, one can compute the Heegaard Floer homology of $HF^+(M,can)$
of $M$  (equivalently, of $-M$) in a purely combinatorial way from
the plumbing graph $\Gammma$. This is true for all the
3-manifolds discussed in this section.  This is done via some
intermediate objects, the {\em graded root} associated with
$\Gammma$ (in fact,  one has a graded root  corresponding to each
$spin^c$-structure of $M$, but here we will discuss only the
`canonical' one). The theory of graded roots, from the point of
view of singularity theory, is rather interesting by itself, and
we plan to exploit further this connection in the future.

Next, we provide a short presentation  of abstract graded roots
(cf. \cite{[61]}).

\subsection{Definition of the `abstract graded root' $(R,\chi)$.} Let $R$ be an
infinite tree with vertices $\calv$ and edges $\cale$. We denote
by $[u,v]$ the edge with
 end-points $u$ and $v$.  We say that $R$ is a graded root
with grading $\chi:\calv\to \Z$ if

(a) $\chi(u)-\chi(v)=\pm 1$ for any $[u,v]\in \cale$;

(b) $\chi(u)>\min\{\chi(v),\chi(w)\}$ for any $[u,v],\
[u,w]\in\cale$;

(c) $\chi$ is bounded below, $\chi^{-1}(n)$ is finite for any
$n\in\Z$, and $\#\chi^{-1}(n)=1$ if $n\gg 0$.

\subsection{Examples.}\label{tau} (1) For any integer $n\in\Z$, let $R_n$ be
the tree with $\calv=\{v^{k} \}_{ k\geq n}$ and
$\cale=\{[v^{k},v^{k+1}]\}_{k\geq n}$. The grading is
$\chi(v^{k})=k$.

(2) Let $I$ be a finite index set. For each $i\in I$ fix  an
integer $n_i\in \Z$; and for each pair $i,j\in I$ fix
$n_{ij}=n_{ji}\in\Z$ with the next properties: (i) $n_{ii}=n_i$;
(ii) $n_{ij}\geq \max\{n_i,n_j\}$; and (iii) $n_{jk}\leq
\max\{n_{ij},n_{ik}\}$ for any $ i,j,k\in I$. For any $i\in I$
consider $R_{n_i}$ with vertices $\{v_i^{k}\}$ and edges
$\{[v_i^{k},v_i^{k+1}]\}$, $(k\geq n_i)$. In the disjoint union
$\coprod_iR_{n_i}$,  for any pair $(i,j)$,
 identify $v_i^{k}$ and $v_j^{k}$,
resp. $[v_i^{k},v_i^{k+1}]$  and $[v_j^{k},v_j^{k+1}]$, whenever
$k\geq n_{ij}$,  and take the induced $\chi$.

(3) Any map $\tau:\{0, 1,\ldots,r\}\to \Z$ produces a starting
data for construction (2). Indeed, set $I=\{0,\ldots,r\}$,
$n_i:=\tau(i)$ ($i\in I$), and $n_{ij}:=\max\{n_k\,:\, i\leq k\leq
j\}$ for $i\leq j$. Then the root constructed in (2) using this
data will be denoted by $(R_\tau,\chi_\tau)$.

\subsubsection{}\label{figure} Here are
two  (typical) graded roots (cf. with \ref{abd}):

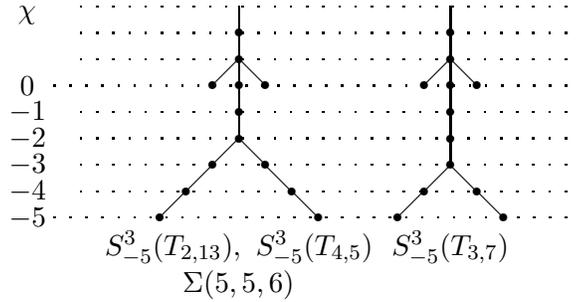
\vspace{3mm}

\begin{picture}(200,130)(-80,290)

\dashline{1}(40,410)(230,410) \dashline{1}(40,400)(230,400)
\put(20,380){\makebox(0,0){0}} \dashline{1}(40,390)(230,390)
\put(20,370){\makebox(0,0){$-1$}} \dashline{1}(30,380)(230,380)
\put(20,360){\makebox(0,0){$-2$}} \dashline{1}(40,370)(230,370)
\put(20,350){\makebox(0,0){$-3$}} \dashline{1}(40,360)(230,360)
\put(20,340){\makebox(0,0){$-4$}} \dashline{1}(40,350)(230,350)
\dashline{1}(40,340)(230,340) \dashline{1}(30,330)(230,330)
\put(20,330){\makebox(0,0){$-5$}} \put(90,380){\circle*{3}}
\put(110,380){\circle*{3}} \put(100,390){\circle*{3}}
\put(100,400){\circle*{3}}
\put(100,360){\circle*{3}}\put(100,370){\circle*{3}}\put(100,380){\circle*{3}}
\put(90,350){\circle*{3}} \put(110,350){\circle*{3}}
\put(80,340){\circle*{3}} \put(120,340){\circle*{3}}
\put(70,330){\circle*{3}} \put(130,330){\circle*{3}}
\put(100,410){\line(0,-1){50}} \put(100,390){\line(1,-1){10}}
\put(100,390){\line(-1,-1){10}} \put(100,360){\line(-1,-1){30}}
\put(100,360){\line(1,-1){30}} \put(170,380){\circle*{3}}
\put(190,380){\circle*{3}} \put(180,390){\circle*{3}}
\put(180,400){\circle*{3}} \put(180,360){\circle*{3}}
\put(180,350){\circle*{3}} \put(180,370){\circle*{3}}
\put(180,380){\circle*{3}} \put(170,340){\circle*{3}}
\put(190,340){\circle*{3}} \put(160,330){\circle*{3}}
\put(200,330){\circle*{3}} \put(180,410){\line(0,-1){60}}
\put(180,390){\line(1,-1){10}} \put(180,390){\line(-1,-1){10}}
\put(180,350){\line(-1,-1){20}} \put(180,350){\line(1,-1){20}}
\put(100,310){\makebox(0,0)[t]{$\Sigma(5,5,6)$}}
\put(180,325){\makebox(0,0)[t]{$S^3_{-5}(T_{3,7})$}}
\put(100,325){\makebox(0,0)[t]{$S^3_{-5}(T_{2,13}),\
S^3_{-5}(T_{4,5})$}} \put(20,410){\makebox(0,0)[t]{$\chi$}}

\end{picture}

\vspace{1mm}

\subsection{The canonical graded root $(R,\chi)$ of $M$.} \cite{[61]} Next, we
define for any (negative definite, plumbed) rational homology
sphere $M$ a graded root.

We fix a plumbing graph $\Gammma$ and denote by $L$ the
corresponding lattice:
 the free $\Z$-module of rank $\#$ with fixed basis
$\{A_j\}_j$, and bilinear form $(A_i,A_j)_{i,j}$. (In our case, a
possible choice is the dual resolution graph and the corresponding
intersection form associated with the minimal good resolution
$Z\to (\{f=0\},0)$.) Set $L'=\Hom_\Z(L,\Z)\subset L\otimes \bq$.
Let $K_Z\in L'$ be the canonical cycle defined by
$K_Z(A_j)+A_j^2+2=0$ for any $j$. Then define $\chi:L\to\Z$ by
(the Riemann-Roch formula) $\chi(x):=-(K_Z(x)+x^2)/2$.

The definition of the graded root captures the  position of the
lattice points in the different ellipsoids $\chi^{-1}(n)$. For any
$n\in \Z$, one constructs a finite 1-dimensional simplicial
complex $\bar{L}_{\leq n}$ as follows. Its 0-skeleton is $L_{\leq
n}:=\{x\in L:\, \chi(x)\leq n\}$. For each $x$ and $j$, with both
$x$ and $x+A_j\in L_{\leq n}$, we consider a unique 1-simplex with
endpoints at $x$ and $x+A_j$ (e.g., the segment $[x,x+A_j]$ in
$L\otimes \R$). We denote the set of connected components of
$\bar{L}_{\leq n}$ by $\pi_0(\bar{L}_{\leq n})$. For any $v\in
\pi_0(\bar{L}_{\leq n})$, let $C_v$ be the corresponding connected
component of $\bar{L}_{\leq n}$.

Next, we define $(R,\chi)$ as follows. The vertices $\calv(R)$ are
$\cup_{n\in \Z} \pi_0(\bar{L}_{\leq n})$. The grading
$\calv(R)\to\Z$, still denoted by $\chi$, is
$\chi|\pi_0(\bar{L}_{\leq n})=n$. If $v_n\in \pi_0(\bar{L}_{\leq
n})$, and $v_{n+1}\in \pi_0(\bar{L}_{\leq n+1})$, and
$C_{v_n}\subset C_{v_{n+1}}$, then $[v_n,v_{n+1}]$ is an edge of
$R$. All the edges are obtained in this way.

\subsection{Example}\label{cusproot} \cite{[63]} Recall that the link of the
superisolated singularity $f$ (where $C$ is rational and
unicuspidal of degree $d$) is the surgery manifold $S^3_{-d}(K)$,
where $K\subset S^3$ is the local knot of $(C,p)$. The graded root
of $M$ can be represented by a function $\tau$ as in \ref{tau}(3)
associated with the Alexander polynomial $\Delta$ of $K\subset
S^3$.  Similarly as in section 6, set $\mu=2\delta$ for the degree
of $\Delta$ (which equals $(d-1)(d-2)$), and write $\Delta(t)$ as
$1+\delta(t-1)+(t-1)^2Q(t)$ for some polynomial $Q(t)=\sum
_{i=0}^{\mu-2} \alpha_it^i$. Set $c_l:=\alpha_{(d-3-l)d}$ (cf. with 6.1). Then
define  $\tau:\{0,1,\ldots, 2d-4\}\to \Z$ \ by

\vspace{2mm}

\hspace{1cm} $\tau(2l) =\frac{l(l-1)}{2}d-l(\delta-1), \ \  \ \
\tau(2l+1)= \tau(2l+2)+c_{d-3-l}.$

\vspace{2mm}

Then $(R,\chi)=(R_\tau,\chi_\tau)$.

\subsection{Example}\label{bri} Let $\Sigma(d,d,d+1)$
be the Seifert 3-manifold $(d,d,d+1)$; equivalently, the link of
the Brieskorn singularity $x^d+y^d+z^{d+1}=0$. Its graded root
also can be represented by the `$\tau$-construction' (for the more
general situation of Seifert manifolds, see \cite{[61]}).

For any $0\leq l\leq d-3$ define $c_l^u:=(l+1)(l+2)/2$, and
$2\delta:=(d-1)(d-2)$. Then define  $\tau^u:\{0,1,\ldots,
2d-4\}\to \Z$ \ by

\vspace{2mm}

\hspace{1cm} $\tau^u(2l) =\frac{l(l-1)}{2}d-l(\delta-1), \ \  \ \
\tau^u(2l+1)= \tau^u(2l+2)+c^u_{d-3-l}.$

\vspace{2mm}

Then $(R,\chi)=(R_{\tau^u},\chi_{\tau^u})$.

\vspace{2mm}

\noindent 
Notice the shocking similarities of \ref{cusproot} and \ref{bri}:
the graded roots associated with $S^3_{-d}(K)$ and $\Sigma(d,d,d+1)$
coincide exactly when $c_l=c^u_l$ for all $l$. 

\vspace{3mm}

\noindent To any graded root, one can associate  a natural graded
$\Z[U]$-module.

\subsection{Definition. The $\Z[U]$-module associated with a graded
root.} \ Consider the $\Z[U]$-module $\Z[U,U^{-1}]$, and
(following \cite{OSzP}) denote by
 $\calt_0^+$ its quotient by the submodule  $U\cdot \Z[U]$.
It is a $\Z[U]$-module with grading  $\deg(U^{-h})=2h$.

Now,  fix a  graded root $(R,\chi)$. Let $\bH(R,\chi)$ be the set
of functions $\phi:\calv\to \calt^+_0$ with the property that
whenever $[v,w]\in \cale$ with $\chi(v)<\chi(w)$ one has $U\cdot
\phi(v)=\phi(w)$. Then $\bH(R,\chi)$ is a $\Z[U]$-module via
$(U\phi)(v)=U\cdot \phi(v)$. Moreover, $\bH(R,\chi)$ has a
grading: $\phi\in \bH(R,\chi)$ is homogeneous of degree $h\in\Z$ if for
each $v\in\calv$ with $\phi(v)\not=0$, $\phi(v)\in\calt^+_0$ is
homogeneous  of degree $h-2\chi(v)$.

\vspace{1mm}

In the sequel, the following notation is useful:  If $P$ is a
$\Q$-graded $\Z[U]$-module with $h$-homogeneous elements $P_h$,
then for any $r\in \Q$  we denote by $P[r]$ the same module graded
in such a way that $P[r]_{h+r}=P_{h}$.

 \subsection{Theorem.} \cite{[61],OSzP}
{\em Assume that $M$ is either $S^3_{-d}(K)$ or $\Sigma(d,d,d+1)$.
Then}
$$HF^+(-M,can)=\bH (R,\chi)[-(K_Z^2+\#)/4].$$

In other words, for these 3-manifolds, the Heegaard-Floer homology
can be recovered from the graded root via $\bH (R,\chi)$ modulo a
shift in grading by $-(K_Z^2+\#)/4$. (The shift in the above two
examples are different; in the case of $S^3_{-d}(K)$ one has
$K_Z^2+\#=1-d(d-2)^2$, while for $\Sigma(d,d,d+1)$ one has
$K_Z^2+\#=-d(d-1)(d-3)$.)

\vspace{2mm}

Now, Conjecture B1 (\ref{CB1}) and the above discussion/examples  read
as follows:

\subsection{Theorem.}\label{last} {\em Assume that $\nu=1$.
Then the following facts are equivalent:

(a) Conjecture B1 (\ref{CB1}) is true,

(b) The canonical graded roots of  $S_{-d}^3(K)$ and
$\Sigma(d,d,d+1)$ are the same.

(c) The canonical Heegaard-Floer homologies of $-S_{-d}^3(K)$ and
$-\Sigma(d,d,d+1)$ are the same modulo a shift in the grading,
namely: }$$HF^+(-S_{-d}^3(K),can)[1-d(d-2)^2]= HF^+(
-\Sigma(d,d,d+1),can)[-d(d-1)(d-3)].$$

\vspace{2mm}

\begin{proof} 
The equivalence $(a) \Leftrightarrow (b)$ is clear from the above discussion,
while  $(b) \Leftrightarrow (c)$ can be deduced by a direct computation, or 
from an easy formula which provides $\bH(R_\tau,\chi_\tau)$ from $\tau$,
cf. \cite{[61]}. (Nevertheless, see another  argument below.)
\end{proof}
\subsection{Remark.} Regarding the Seiberg-Witten invariant of
$M=S^3_{-d}(K)$, one has
\begin{equation*}\ssw_M(can)-\frac{K_Z^2+\#}{8}=\sum_{l\geq
0}\tau(2l+1)-\tau(2l+2)=\sum_{l\geq 0}c_l,\tag{29}\end{equation*}
 and there is a
similar formula for $M=\Sigma(d,d,d+1)$ with the obvious
replacements.

Therefore, to the equivalences (\ref{last}) one can add:
 $$(d)\ \ \
\ssw_M(can)-\frac{K_Z^2+\#}{8}\
|_{M=S^3_{-d}(K)}=\ssw_M(can)-\frac{K_Z^2+\#}{8}\
|_{M=\Sigma(d,d,d+1)}.$$ Since the Conjecture $SWC$ (\ref{SWC}) is
true for the Brieskorn singularity $f_{BR}:=x^d+y^d+z^{d+1}$ (cf.
\cite{[52]}), and the geometric genus of the superisolated
singularity $f$ equals the geometric genus of $f_{BR}$ (both equal
$d(d-1)(d-2)/6$),  this last  identity $(d)$ is also equivalent
with the validity of the $SWC$ for $f$ --- \, a fact already
proved in \ref{top}.

(Notice also that the expression 
$\ssw_M(can)-(K_Z^2+\#)/8$ can be deduced from $\bH$, a fact which
implies  $(c)\Rightarrow 
(d)$,  while  $(d)\Rightarrow (a)$ follows from  (29).)  

\subsection{Example.}\label{abd} Assume that $d=5$ and $C$ is unicuspidal
whose singular point has only one Puiseux pair $(a,b)$ with $a<b$.
Then by the genus formula the possible values of $(a,b)$ are
$(4,5)$, $(3,7)$ and $(2,13)$. It turns out that the first and the
third cases can be realized, while the second not. The
corresponding graded roots (together with the root of
$\Sigma(5,5,6)$) are drawn in the above figure (\ref{figure}).

\end{document}